\documentclass[english]{smfart} 
\usepackage{smfthm} 

\DeclareMathOperator{\ad}{ad}

\DeclareMathOperator{\End}{End}
\DeclareMathOperator{\Hom}{Hom}

\newcommand{\la}{\lambda}
\newcommand{\ot}{\otimes}
\newcommand{\proj}{\mathrm proj}

\newcommand{\ra}{\rightarrow}
\newcommand{\sig}{\sigma}
\newcommand{\sq}[1]{#1^\times/#1^{\times 2}}

\begin{document}
\frontmatter
\title[A generic characterization of direct summands] {A generic characterization of direct summands for 
  orthogonal involutions}
\author{Anne Qu\'eguiner-Mathieu}
\address{Laboratoire Analyse, G\'eom\'etrie \& Applications\\
UMR CNRS 7539 -- Institut Galil\'ee\\
Universit\'e Paris 13\\
93430 Villetaneuse\\
France}
\email{queguin@math.univ-paris13.fr}
\urladdr{http://www-math.univ-paris13.fr/{\textasciitilde}queguin/}
\date{\today}
\begin{abstract}
The `transcendental methods' in the algebraic theory of quadratic
forms are based on two major results, proved in the 60's by Cassels
and Pfister, and known as the representation and the subform theorems. 
A generalization of the 
representation theorem was
proven by Jean-Pierre Tignol in 1996, in the setting of central simple
algebras with involution. This paper studies the subform
question for orthogonal involutions. A generic
characterization of direct summands is given; an analogue of the
subform theorem is proven for division algebras and algebras of index
at most $2$. 
\end{abstract}
\maketitle
\mainmatter
\section*{Introduction} 
The `transcendental methods' in the algebraic theory of quadratic
forms are based on two major results, proved in the 60's by Cassels
and Pfister, and known as the representation and the subform theorems
(see~\cite[Ch.~4 \S3]{Sch} or \cite[CH.~9 \S1 and 2]{Lam2}). 
J.-P. Tignol gave in~\cite{T:96} a generalization of the 
representation theorem for algebras with involution (of any kind), 
which implies the corresponding statement for quadratic forms. 

In this paper, the subform question is studied in the context of
algebras with orthogonal involutions. The main results, which 
give partial answers to this question, are 
theorems~\ref{scai.theo} and~\ref{scaibis.theo}, stated and proved in
\S\ref{scai.sec} and~\S\ref{scaibis.sec}. 
The first one gives a generic characterization of direct summands,
which is valid for any algebra with orthogonal involution, but which is much
weaker than the subform theorem in the split case. 
The second one is an analogue of the subform theorem, but only for division
algebras and algebras of index at most $2$. 
Before proving these theorems, we define in~\S\ref{directsummands.sec}
direct summands of an algebra with involution, 
using the direct sum of~\cite{D}, and we restate the subform
theorem in a convenient way for our purpose
in~\S\ref{subform.sec}. 

We assume throughout the paper that the base field $F$ has
characteristic different from~$2$, and refer the reader to~\cite{KMRT}
for basic facts on algebras with involution. 
\section{Direct summands of an algebra with involution}
\label{directsummands.sec} 
Consider two central simple algebras with involution $(B,\tau)$
and $(B',\tau')$ over $F$, which 
are Morita equivalent, i.e. $B$ and $B'$ are Brauer equivalent and $\tau$ 
and $\tau'$ 
are of the same type. 
They can be represented as $(B,\tau)=(\End_D(N),\ad_{h_N})$ and 
$(B',\tau')=(\End_D(N'),\ad_{h_{N'}})$ for some hermitian modules $(N,h_N)$
and $(N',h_{N'})$ over a central division algebra with involution $(D,\gamma)$ over $F$. 
The direct sum for hermitian modules gives rise to an algebra with
involution, $(\End_D(N\oplus N'),\ad_{h_N\oplus h_{N'}})$,
which is Morita equivalent to
$(B,\tau)$ and $(B',\tau')$. 
But the involutions $\tau$ and $\tau'$ only determine the hermitian
forms $h_N$ and $h_{N'}$ up to a scalar factor, and modifying this factors
independently may lead to a hermitian form which is not similar to
$h_N\oplus h_{N'}$. 

In~\cite{D}, Dejaiffe defined a notion of direct sum for algebras with
involution which extends direct sum of hermitian modules, i.e. such
that $(\End_D(N\oplus N'),\ad_{h_N\oplus h_{N'}})$ is a direct sum of
$(\End_D(N),\ad_{h_N})$ and $(\End_D(N'),\ad_{h_{N'}})$. 
Precisely, given any Morita equivalence data between two
algebras with involution $(B,\tau)$ and $(B',\tau')$ she defines the
corresponding direct sum of $(B,\tau)$ and $(B',\tau')$. Note that different Morita equivalence datas
(which amount to modifying scalars as in the previous paragraph) may lead
to non isomorphic direct sums of the same $(B,\tau)$ and $(B',\tau')$. 

We say that $(B,\tau)$ is a {\em direct summand} of
$(A,\sigma)$ if there exist $(B',\tau')$, Morita equivalent to
$(B,\tau)$ and a direct sum of $(B,\tau)$ and $(B',\tau')$ which is
isomorphic to $(A,\sigma)$. 
By~\cite[Prop.~2.2]{D}, this leads to the following definition:  
\begin{defi} The algebra with involution $(B,\tau)$ is a direct summand of $(A,\sigma)$   
if there exist a $\sigma$-symmetric idempotent $e\in A$ such that 
$(eAe,\sigma_{|eAe})$ is isomorphic to $(B,\tau)$. 
\end{defi}
This condition can be translated in terms of hermitian modules as
follows. Fix representations $(A,\sigma)=(\End_D(M),\ad_{h_M})$ and
$(B,\tau)=(\End_D(N),\ad_{h_N})$, for some hermitian modules over a central division algebra
with orthogonal involution $(D,\gamma)$ over $F$. 
Denote by $M_0\subset M$ the image of the idempotent $e$, and by
$h_{M_0}$ the restriction of $h_M$ to $M_0$. 
The algebra with involution $(eAe,\sigma_{|eAe})$ is
$(\End_D(M_0),\ad_{h_{M_0}})$. 
Hence, $(B,\tau)$ is a
direct summand of $(A,\sigma)$ if and only if 
$M$ contains a submodule $M_0$ such that $h_{M_0}$ is similar to $h_M$. 

In the split orthogonal case, 
that is $(A,\sigma)=(\End_F(V),\ad_{q_V})$ and
$(B,\tau)=(\End_F(W),\ad_{q_W})$, 
for some quadratic spaces $(V,q_V)$ and $(W,q_W)$ over $F$, 
we get that $(B,\tau)$ is a direct summand of $(A,\sigma)$ 
if and only if there exists a scalar $\la\in F^\times$ such that 
$\la q_W$ is a subform of $q_V$. 
A 'generic' condition under which this is satisfied is given by the
subform theorem, at least in a version up to similarities. 
For further use, we give in the next section a projective version of this statement. 
\section{A projective version of the subform theorem up to similarities}
\label{subform.sec} 
\subsection{Definition of ${q_W}^{proj}$}
\label{{q_W}proj.subsection}
Denote by $F(W)$ the function 
field of the affine variety 
$W$ and by $F({\mathbb P}W)$ the field of rational functions on the 
projective variety ${\mathbb P}W$. 
The generic point of the projective space ${\mathbb P}W$, viewed as an
$F({\mathbb P}W)$ rational point, gives rise to a line $L_W\subset
W_{F({\mathbb P}W)}$, which we call the generic line. 
We define the {\em projective class of $q_W$} to be the square class in
$F({\mathbb P}W)$ of the value of $q_W$ at any point of the generic
line. 
It is a well defined element of $\sq {F({\mathbb P}W)}$, and the
notation ${q_W}^{\proj}$ stands for an element in $F({\mathbb P}W)$ who 
belongs to the projective class of ${q_W}$.

If we identify $F({\mathbb P}W)$ with 
the subset $F(W)_0\subset F(W)$ of 
degree $0$ homogeneous functions, we may describe
the projective class of $q_W$ as the quotient 
$\frac{{q_W}}{f^2}$, where $q_W$ is viewed as an element of $F(W)$ and
$f$ is an arbitrary degree $1$ homogeneous element so that the
quotient is in $F(W)_0$. The square class of $\frac{{q_W}}{f^2}$
clearly does not depend on the choice of $f$. 

Let $(e_1,\dots,e_n)$ be any basis of $W$ over $F$, with dual basis 
$t_1,\dots t_n$, so that $F(W)\simeq F(t_1,\dots,t_n)$. 
We may then identify $F({\mathbb P}W)$ with $F(\theta_2,\dots
\theta_n)$, where
$\theta_2,\dots,\theta_n$ are indeterminates, by  
$f\mapsto f(1,\theta_2,\dots,\theta_n)$. 
Assume moreover that the basis $(e_1,\dots e_n)$ is 
orthogonal for ${q_W}$, and let $b_i={q_W}(e_i)$. 
Using these identifications, one may check that the element $b_1+b_2\theta_2^2+\dots+b_n\theta_n^2$
belongs to the projective class of ${q_W}$. 
Indeed, we have 
$b_1+b_2\theta_2^2+\dots+b_n\theta_n^2={q_W}(\delta)$, for
$\delta=e_1+e_2\theta_2+\dots+e_n\theta_n\in L_W\subset W_{F({\mathbb P}W)}$. 
\subsection{A projective version of the subform theorem up to similarities} 
The classical subform theorem (\cite[Ch.~4, Th.~3.7]{Sch},
\cite[Ch.~9, Th.~2.8]{Lam2}) has the following easy consequence: 
\begin{prop} 
\label{projectivesubform.theo}
Let $(V,q_V)$ and $(W,{q_W})$ be two quadratic spaces, with $q_V$ anisotropic. 
There exists $\la\in F^\times$ such that $\la {q_W}$ is a subform of $q_V$ 
if and only if 
\begin{equation}\label{dep}\tag{dep}
\hspace{2cm}\exists\la\in F^\times
\mbox{ such that }{q_V}_{F({\mathbb P}W)}\mbox{ represents }\la {q_W}^{\proj}. 
\end{equation}
\end{prop} 
\begin{rema} 
This statement does not depend on the choice of ${q_W}^{proj}$, since the 
set of values represented by a quadratic form is stable under multiplication 
by a square.  
\end{rema} 
\begin{proof} 
Since ${q_W}_{F({\mathbb P}W)}$ represents $q_W^{\proj}$ (see
\S\ref{{q_W}proj.subsection}), 
condition~(\ref{dep}) clearly is necessary. 
Let us prove it is also
sufficient. 
Pick a vector 
$v\in V_{F({\mathbb P}W)}$ satisfying 
${q_V}_{F({\mathbb P}W)}(v)=\la(b_1+b_2\theta_2^2+\dots+b_n\theta_n^2)$. 
If we identify $F({\mathbb P}W)(t)=F(\theta_2,\dots,\theta_n)(t)$ with $F(W)$ by 
$f\mapsto f(\frac{t_2}{t_1},\dots,\frac{t_n}{t_1})(t_1)$, we get 
${q_V}_{F(W)}(tv)=\la(b_1t_1^2+\dots+b_nt_n^2)$. 
By the classical subform theorem, this implies that $\la {q_W}$ is a subform of 
$q_V$. 
\end{proof} 
\begin{rema}
\label{projectivesubform.rema} 
There are two differences between this statement and the classical 
subform theorem. 
Using $F({\mathbb P}W)$ instead of $F(W)$ and $q_W^{proj}$ instead of $q_W$ 
is not a serious 
one. 
As opposed to this, one should notice that it does not seem easy to prove the 
subform theorem 
up to similarities (even in an affine version) without using the subform theorem. 
\end{rema} 

\section{The generic ideal and a characterization of direct summands} 
\label{scai.sec}

We now go back to the setting of algebras with involution, and for
simplicity, we restrict ourselves to the orthogonal case. 
Hence $(A,\sigma)$ and $(B,\tau)$ are Brauer equivalent central simple
algebras over $F$, with orthogonal involutions. 

We denote by $X_B$ the Brauer-Severi variety of $B$ 
and by $F_B$ the field of rational functions on $X_B$. 
It is a generic splitting field 
for the algebra $B$. 
By definition of $X_B$, its generic point, 
viewed as an $F_B$-rational point, corresponds to 
a right ideal of reduced dimension $1$ 
of the split algebra $B_{F_B}:=B\otimes_F F_B$. 
We call it 
the {\em generic ideal} of $B$ and denote it by $I_B$. 
By~\cite[(1.12)]{KMRT}, given a representation $B_{F_B}=\End_{F_B}(W)$, there exists 
a unique line $L_B$ in the $F_B$-vector space $W$
such that $I_B=\Hom_{F_B}(W,L_B)$; we call it the {\em generic line}
of $B$. 

If $B$ is split and $(B,\tau)=(\End_F(W),\ad_{q_W})$ for some
quadratic space $(W,q_W)$ over $F$, then $X_B$ is isomorphic to 
${\mathbb P}W$, and with the notations of \S{\ref{{q_W}proj.subsection}}, we
have $L_B=L_W$ and 
$I_B=Hom_{F({\mathbb P}W)}(W,L_W)$. 

In this section, we prove the following theorem, which gives 
a necessary and sufficient condition under which 
$(A,\sigma)$ contains $(B,\tau)$ as a direct summand in terms of the generic ideal 
of $B$. Note that we do not view this as 
an analogue of the subform theorem (see remark~\ref{pasanalogue.rema} 
below). In particular, the involution $\sigma$ may be isotropic. 

\begin{theo} \label{scai.theo}
Let $(A,\sigma)$ and $(B,\tau)$ be two Brauer-equivalent central 
simple algebras, 
each endowed with an orthogonal involution. 
Then $(A,\sigma)$ contains $(B,\tau)$ as a direct summand if and only
if 
\begin{multline} \label{algebra.cond}\tag{iso}
\mbox{there exist a }\sigma\mbox{-symmetric idempotent }e\in
A \\
\mbox{ and an isomorphism }\Psi:\,eAe\tilde\ra B
\mbox{ such that }\\
\forall g\in I_B,\ \exists f\in A_{F_B},\ 
\Psi_{F_B}(e\sig(f)fe)=\tau(g) g.\\
\end{multline}
\end{theo} 

\begin{rema} (i) The idea of considering elements of the type $\tau(g)g$,
  with $g$ of rank $1$, is borrowed from Tignol's paper~\cite{T:96}. 
Since $B_{F_B}$ is split, such an element may give us some information on the
value of the underlying $F_B$-quadratic form on the image of $g$,
namely the generic line (see lemma~\ref{tau(g)g.lemm} below). In the
split case, this value precisely is the projective class of $q_W$. 

(ii) Since $A$ and $B$ are Brauer equivalent, there exist in general 
many $\sigma$-symmetric idempotents $e\in A$ such that $eAe$ is
isomorphic to $B$. Indeed, one may take for $e$ any orthogonal projection on a submodule of
$M$ of dimension over $D$ equal to $\dim_D(N)$. 
Given such an idempotent $e$, Theorem~\ref{scai.theo} actually gives a
criterion of isomorphism between the involutions $\sig_{|eAe}$ and
$\tau$.  
\end{rema} 

\begin{proof}
Condition~(\ref{algebra.cond}) is clearly necessary. 
Indeed, if $(B,\tau)$ is a direct summand of $(A,\sigma)$, 
there exist $e\in A$ such that $e^2=\sig(e)=e$ and 
an isomorphism of algebras with involution 
$\Psi:\,(eAe,\sig_{|eAe})\tilde\ra (B,\tau)$. 
For any $g\in I_B\subset B_{F_B}$, the element 
$f=\psi_{F_B}^{-1}(g)$ clearly satisfies the required condition. 

We have to prove condition~(\ref{algebra.cond}) is also sufficient.  
Consider $(A,\sigma)$ and $(B,\tau)$ as in the theorem, and take
representations $(A,\sigma)=(\End_D(M),\ad_{h_M})$ and
$(B,\tau)=(\End_D(N),\ad_{h_N})$, for some hermitian modules over a central division algebra
with orthogonal involution $(D,\gamma)$ over $F$. 
Let $M_0\subset M$ be the image of $e$, so that $e$ is the orthogonal projection 
on $M_0$. The isomorphim $eAe=\End_D(M_0)\simeq B$ is given by some 
isomorphism of $D$-modules $\psi:\,M_0\tilde\ra N$. 
Let us denote by $h_{M_0}$ the restriction of $h_M$ to $M_0$. 
We will prove that $h_{M_0}$ is similar to $h_N$, which in turn 
implies that $(B,\tau)$ is a direct summand of $(A,\sigma)$. 

Let $n_1,\dots,n_s$ be an orthogonal basis of $(N,h_N)$. 
The elements $m_k=\psi^{-1}(n_k)$ for $1\leq k\leq s$ 
form a basis of the $D$-module $M_0$. 
We will actually compute $h_{M_0}$ in this
  basis, using condition~(\ref{algebra.cond}).

Since $A$ and $D$ are Brauer equivalent to $B$, they both split over $F_B$. 
Let $(E,q_{\gamma})$ be a quadratic space over $F_B$ such that 
$(D_{F_B},\gamma)\simeq(\End_{F_B}(E),\ad_{q_{\gamma}})$. 
By Morita theory, we then have (see for instance~\cite[\S1.4]{BP})
\[B_{F_B}\simeq \End_{F_B}(W)\text{, where }
W=(N\otimes_F F_B)\otimes_{D\otimes_F F_B} E,
\] and similarly, 
\[A_{F_B}\simeq \End_{F_B}(V)\text{, 
where }V=(M\otimes_F F_B)\otimes_{D\otimes_F F_B} E\text{, and }
\] 
\[eA_{F_B}e\simeq \End_{F_B}(V_0)\text{, 
where }V_0=(M_0\otimes_F F_B)\otimes_{D\otimes_F F_B} E.\] 
Moreover, the involutions $\sigma$, 
$\tau$ and $\sigma_{|eAe}$ are respectively adjoint, after 
scalar extension to $F_B$ to the quadratic forms $q_V$, $q_W$ and $q_{V_0}$ 
defined by 
\[b_{q_V}((m\otimes\lambda)\otimes e,(m'\otimes\lambda')\otimes e')=
b_{q_\gamma}(e,(h_M(m,m')\otimes\la\la')(e')),\] 
and similarly for $q_W$ and $q_{V_0}$.

Let $e_1,\dots,e_d$ be an orthogonal basis of $(E,q_\gamma)$. 
We denote by $d_i=h_N(n_i,n_i)$ and $a_j=q_\gamma(e_j)$. The elements 
$n_i\otimes1\otimes e_j$ for $1\leq i\leq s$ and $1\leq j\leq d$ form 
a basis of the $F_B$ vector space $W$. 
Let $\delta$ be any non trivial element in the 
generic line $L_B\subset W$ defined at the beginning of this section. 

The proof is based on the following computation, which the reader may
easily check: 
\begin{lemm} 
\label{tau(g)g.lemm}
Consider the element $g\in I_B$ defined by 
$g(n_i\otimes 1\otimes e_j)=\alpha_{i,j}\delta$. 
Then, $\tau(g)g$ maps $n_i\otimes 1\otimes e_j$ to 
$\alpha_{i,j}\,q_W(\delta)(\sum_{1\leq k\leq s\,,\,1\leq l\leq d} 
\alpha_{k,l}(n_kd_k^{-1}\otimes 
1\otimes\frac{e_l}{a_l}))$. 
\end{lemm}

For any couple $1\leq i\leq s$ and 
$1\leq j\leq d$, let us first apply this lemma to the element $g_{i,j}\in I_B$ 
which maps 
$n_i\otimes 1\otimes e_j$ to $\delta$, and any other element of 
the basis to $0$, so that 
$\tau(g_{i,j})g_{i,j}$ maps $n_i\otimes 1\otimes e_j$ to 
$n_id_i^{-1}\otimes q_W(\delta)\otimes \frac{e_j}{a_j}$ and any other element 
of the basis to $0$. 
Denote by $f_{i,j}$ the corresponding element of $A_{F_B}=\End_{F_B}(V)$ 
given by condition~(\ref{algebra.cond}). It satisfies 
\[\psi_{F_B}e\sig(f_{i,j})f_{i,j}e\psi_{F_B}^{-1}=\tau(g_{i,j})g_{i,j}.\]
from which we deduce that $e\sigma(f_{i,j})f_{i,j}e$ maps 
$m_i\otimes 1\otimes e_j$ to $m_id_i^{-1}\otimes q_W(\delta)\otimes \frac{e_j}{a_j}$ 
and 
any other element of the basis to $0$.  
Now, we can compute 
\[b_{q_{V_0}}(f_{i,j}(m_k\otimes 1\otimes e_l),f_{i,j}(m_p\otimes 1\otimes e_q))\] 
for any $k,p\in\{1,\dots,s\}$ and $l,q\in\{1,\dots,d\}$ in two different ways. 
First, it is equal to 

\begin{multline} 
\label{q(nu).equa}
b_{q_{V_0}}(m_k\otimes 1\otimes
e_l,e\sigma(f_{i,j})f_{i,j}e(m_p\otimes 1\otimes e_q)) \\
=\left\{
\begin{array}{l}
0\text{ if }p\not=i\text{ or }q\not=j\\
q_W(\delta)b_{q_\gamma}(e_l,(h_{M_0}(m_k,m_i)d_i^{-1})(\frac{e_j}{a_j}))
\text{ if }(p,q)=(i,j).\\
\end{array}\right. \\
\end{multline} 

By symmetry, it is also equal to 
$\left\{\begin{array}{l}
0\text{ if }k\not=i\text{ or }l\not=j\\
q_W(\delta)b_{q_\gamma}(e_q,(h_{M_0}(m_p,m_i)d_i^{-1})(\frac{e_j}{a_j}))
\text{ if }(k,l)=(i,j).\\
\end{array}\right.$

From this, we deduce that if $k\not=i$, then for any $l$, 
\[b_{q_\gamma}(e_l,(h_{M_0}(m_k,m_i)d_i^{-1})(e_j))=0.\]
Hence $(h_{M_0}(m_k,m_i)d_i^{-1})(e_j)=0$. This is valid for any value of $j$, 
and we finally get 
\begin{equation}
\label{orthogonal.equa}
h_{M_0}(m_k,m_i)=0\text{ if }k\not=i
\end{equation}

Let us take now $k=i$. 
For any $l\not=j$, we have 
\[b_{q_\gamma}(e_l,(h_{M_0}(m_i,m_i)d_i^{-1})(e_j))=0.\]
Since the basis $(e_1,\dots,e_d)$ is orthogonal for $q_\gamma$, 
this implies that there exists an element $\lambda_{i,j}\in F_B^\times$ 
such that 
\[(h_{M_0}(m_i,m_i)d_i^{-1})(e_j))=\lambda_{i,j} e_j.\]
In other words, the element 
$h_{M_0}(m_i,m_i)d_i^{-1}\in D\subset D\otimes_F F_B\simeq\End_{F_B}(E)$ 
corresponds in the basis $e_1,\dots,e_d$ to the diagonal matrix with coefficients 
$\lambda_{i,j}$, $1\leq j\leq d$. 

Let us now prove that the coefficients $\lambda_{i,j}$ are all equal. 
Consider the element 
$g_i\in I_B$ which maps $n_i\otimes 1\otimes e_j$ to $\delta$ 
for all $j=1,\dots,d$ and any other element of the basis to $0$. 
Again by the previous lemma, $\tau(g_i)g_i$ maps $n_i\otimes 1\otimes e_j$ to 
$n_id_i^{-1}\otimes q_W(\delta)\otimes
(\frac{e_1}{a_1}+\dots+\frac{e_d}{a_d})$ for any $j$, 
and any other element of the basis to $0$.  
Let $f_i$ be the corresponding element of $A_{F_B}\simeq\End_{F_B}(V)$ 
given by condition~(\ref{algebra.cond}). 
We get that $e\sig(f_i)f_ie$ maps any $m_i\otimes 1\otimes e_j$ to 
$m_id_i^{-1}\otimes q_W(\delta)\ot (\frac{e_1}{a_1}+\dots+\frac{e_d}{a_d})$.
The same computation as above for 
\[b_{q_{V_0}}(f_i(m_i\otimes 1\otimes e_j),f_i(m_i\otimes 1\otimes e_l)),\] 
now gives $\lambda_{i,j}=\lambda_{i,l}$, which proves that 
$h_{M_0}(m_i,m_i)d_i^{-1}$ is actually central, $h_{M_0}(m_i,m_i)d_i^{-1}=
\lambda_i\in F_B^\times$. 
Since we know from the very beginning it lies in $D\subset D\otimes_F F_B$,  
$\lambda_i$ actually belongs to $F^\times$. 

To finish the proof, consider the element $g\in I_B$ which maps any 
element of the basis to $\delta$. 
Then, $\tau(g)g$ maps any $n_i\otimes 1\otimes e_j$ to 
\[(n_1d_1^{-1}+\dots+n_{s}d_{s}^{-1})\otimes q_W(\delta)\otimes 
(\frac{e_1}{a_1}+\dots+\frac{e_d}{a_d}).\]
Again, let $f$ be the corresponding element of $A_{F_B}\simeq\End_{F_B}(V)$ 
given by condition~(\ref{algebra.cond}). 
We get that $e\sigma(f)fe$ maps any $m_i\otimes 1\otimes e_j$ to 
\[(m_1d_1^{-1}+\dots+m_{s}d_{s}^{-1})\otimes q_W(\delta)\otimes 
(\frac{e_1}{a_1}+\dots+\frac{e_d}{a_d}).\]
Computing as before $b_q(f(m_i\otimes 1\otimes e_j),f(m_k\otimes 1\otimes e_l))$ 
in two different ways, we get 
$\lambda_i=\lambda_k=\lambda\in F^\times$. 
Hence we have proven that 
\begin{equation}
\label{similarity.equa}
h_{M_0}(m_i,m_i)=\lambda d_i=\lambda h_N(n_i,n_i)\text{ for any }1\leq i\leq s,
\end{equation}
and combined with~\ref{orthogonal.equa}, this finishes the proof. 
\end{proof} 

\section{A subform theorem in some particular cases} 
\label{scaibis.sec}

Under some assumption on the algebra, we may improve the
previous theorem. Note that conditions (ii) and (iii) in the following
statement are direct consequences of condition~(\ref{algebra.cond}) of
theorem~\ref{scai.theo}. Also, the involution $\sigma$ now is supposed
to be anisotropic.

\begin{theo} \label{scaibis.theo}
Let $(A,\sigma)$ and $(B,\tau)$ be two Brauer-equivalent central 
simple algebras, 
each endowed with an orthogonal involution, with $(A,\sigma)$ anisotropic. 
We assume moreover that either the index $d$ of the algebras $A$ and $B$ is 
at most $2$ or $B$ is a division algebra. 
Then the following assertions are equivalent : 
\begin{multline*}\tag{i}
(B,\tau)\mbox{ is a direct summand of }(A,\sigma); 
\end{multline*} 
\begin{multline*} \tag{ii}
\mbox{For any left ideal }J\subset B\mbox{ of reduced dimension }d
\mbox{ with }\tau(J) 
\mbox{ anisotropic,}\\
\mbox{ there exists }
\text{a symmetric idempotent }e\in A
\mbox{ and an isomorphism }\Psi: eAe\simeq B\\
\mbox{ such that }
\forall g\in I_B\cap J_{F_B},\ \exists f\in A_{F_B},\ \text{with }
\Psi_{F_B}(e\sig(f)fe)=\tau(g)g;
\end{multline*} 
\begin{multline*} \tag{iii}
\text{There exists a left ideal }J\subset B\text{ of reduced dimension }d
\text{ with }\tau(J)
\text{ anisotropic,}\\
\text{a symmetric idempotent }e\in A\text{, and an isomorphism }\Psi: eAe\simeq B,\\
\text{such that }
\forall g\in I_B\cap J_{F_B},\ \exists f\in A_{F_B},\ \text{with }
\Psi_{F_B}(e\sig(f)fe)=\tau(g)g.
\end{multline*} 
\end{theo}

\begin{rema}
\label{ideal.rema} 
It follows from~\cite[(1.12) and (6.2)]{KMRT} that the algebra $B$ 
always contains left ideals $J$ of reduced 
dimension $d$ such that $\tau(J)$ is anisotropic. 
Indeed, $J=\Hom_D(N/{(n.D)^\perp},N)$ satisfies these conditions as
soon as the vector $n\in N$ is anisotropic. 

Moreover, for any such ideal $J$, after scalar extension $J_{F_B}$ coincides 
with the set of endomorphisms 
of the split algebra $B_{F_B}=\End_{F_B}(W)$ whose kernel contains 
some particular subspace 
$W_1\subset W$ of codimension $d$. 
Hence, since $d$ is non zero and $L_B$ has dimension $1$, there exists non trivial 
elements in $J_{F_B}$ with image included in $L_B$, that is non trivial 
elements in the 
intersection $J_{F_B}\cap I_B$. 
This guarantees conditions (ii) and (iii) are non empty. 
\end{rema}

Clearly, (i) implies (ii) and (ii) implies (iii). 
We now prove the theorem in the split and division cases. 
The index $2$ case will be proven in \S~\ref{ind2.sect}. 

\subsection{The division case} 

In the particular case when $B$ is a division algebra, there exists 
a unique left ideal 
of reduced dimension 
$d$, $B$ itself, which is necessarily anisotropic. 
Hence conditions (ii) and (iii) in that case both reduce to 
condition~(\ref{algebra.cond}) of theorem~\ref{scai.theo}. 
This already proves the result is true if 
$B$ is a division algebra. 
Note that the anisotropy hypothesis is not necessary 
in that case.

\subsection{The split case} 

We assume now that $A$ and $B$ are split, and take representations 
$(A,\sigma)=(\End_F(V),\ad{q_V})$ and $(B,\tau)=(\End_F(W),\ad_{q_W})$. 
The involution $\sigma$ is anisotropic if and only if ${q_V}$ is anisotropic, and 
$(A,\sigma)$ contains $(B,\tau)$ as a direct summand if and only if ${q_V}$ contains 
$\la q_W$ as a subform for some $\la\in F^\times$. 

Hence the theorem follows from the following proposition : 
\begin{prop} 
\label{split.prop}
If $A$ and $B$ are split, then condition (iii) is equivalent to
condition~(\ref{dep}) of Theorem~\ref{projectivesubform.theo}. 
\end{prop} 

\begin{proof}
Recall the isomorphism $B=\End_F(W)\simeq W\otimes W$ given by 
$(x\otimes y)(z)=x b_{q_W}(y,z)$. 
It is an isomorphism of algebras with involution if we endow $W\otimes W$ 
with the product $(x\otimes y)(z\otimes w)=b_{q_W}(y,z) x\otimes w$ and 
the involution $\tau(x\otimes y)=y\otimes x$. 

Under this isomorphism, a left ideal $J\subset B=\End_F(W)$ of reduced dimension $1$ 
corresponds to $\{x\otimes w,\ x\in W\}$ for some non trivial $w\in W$, 
uniquely defined up to a scalar factor. 
Moreover, $\sig(J)=\{w\otimes x,\ x\in W\}$ is anisotropic if and only if 
the vector $w$ itself is anisotropic. 
On the other hand, after scalar extension to $F_B$, 
the ideal $I_B$ corresponds under the same isomorphism to $L_B\otimes _{F_B}W_{F_B}$. 

Let us first assume (iii).  
Denote by $V_0$ the image of the idempotent $e$. 
The isomorphims $\Psi:\,eAe=\End_F(V_0)\ra B$ is given by
$\Psi(f)=\psi f \psi^{-1}$, 
for some isomorphism $\psi:\,V_0\simeq W$. 

Fix a vector $w\in W$ such that $J=\{x\otimes w,\ x\in W\}$. 
Any $g\in I_B\bigcap J_{F_B}$ can be written as 
$g=\delta\otimes w$ for some $\delta\in L_B$.  
We then have 
$\tau(g)g=q_W(\delta)w\otimes w$. 
Hence, $\tau(g)g$ maps $w$ to $q_W(\delta)q_W(w)w$, and 
any element of the orthogonal of $w$ to $0$. 

Let us now consider the element $v=\psi^{-1}(w)\in V_0$, and denote by $f$ 
the element of $\End_{F_B}(V_{F_B})$ corresponding to $g$ given by condition (iii). 
We have 

${q_V}(f(v))={q_V}(fe(v))=b_{q_V}(e\sigma(f)fe\psi^{-1}(w),v)=b_{q_V}(\psi^{-1}
(\tau(g)g(w)),v)
={q}_W(\delta)q_W(w){q_V}(v)$. 
This proves (\ref{dep}) is satisfied. Indeed, we have already noticed in
\S\ref{scai.sec} that $q_W(\delta)$ 
belongs to the projective class of $q_W$. On the other hand, $v$ and $w$ are both 
defined over $F$, $w$ is anisotropic and 
${q_V}$ is anisotropic. 
Hence $q_W(w){q_V}(v)\in F^\times$. 

\smallskip

Let us assume conversely that (\ref{dep}) is satisfied, and let
$\nu\in V_{F_B}$ be a vector satisfying 
${q_V}(\nu)=\la q_W^{proj}=\la q_W(\delta_0)$ for some fixed non trivial 
$\delta_0\in L$. 
Specializing this equality, 
one may find anisotropic vectors 
$w\in W$ and $v\in V$ such that $q_V(v)=\la q_W(w)$. 
We let $v_1=v,v_2,\dots,v_m$ be an orthogonal basis of $(V,{q_V})$ 
and $w_1=w,w_2,\dots,w_n$ be an orthogonal basis of $(W,q_W)$. 
Let $V_0\subset V$ be the sub-vector space generated by $v_1,\dots v_n$.  
We define $e$ to be the orthogonal projection on $V_0$, 
$\psi:\,V_0\tilde\ra W$ the isomorphism defined by $\psi(v_i)=w_i$ for $i=1,\dots,n$, 
and $\Psi:\,eAe\tilde\ra B$ the isomorphism defined by 
$\Psi(f)=\psi f\psi^{-1}$.  

Any element $g\in I_B\cap J_{F_B}$ can be written as 
$g=l\delta_0\otimes w$ for some $l\in F_B$. 
Take the corresponding element in $A_{F_B}$ to be $f=\frac{l}{\la}\nu\otimes v$. 
We then have 
$\tau(g)g=l^{2}q_W(\delta_0)w\otimes w$ 
while $\sigma(f)f=\frac{l^2}{\la^2}{q_V}(\nu)v\otimes v$. 
An easy computation shows that 
$\psi e\sigma(f)fe\psi^{-1}(w)=\frac{l^2}{\la^2}{q_V}(\nu){q_V}(v) w=
l^2 q_W(\delta_0)q_W(w)w=\tau(g)g(w)$. 
Moreover, both maps are trivial on the orthogonal of $w$, 
since by definition of $\psi$, the image under $\psi^{-1}$ of 
the orthogonal of $w$ 
is included in the orthogonal of $v$. 
Hence $\psi e\sigma(f)fe\psi^{-1}=\tau(g)g$, and this proves (iii).  
\end{proof}

\begin{rema} 
\label{pasanalogue.rema}
Note that proposition~\ref{split.prop} actually shows that, in the split case, 
theorem~\ref{scaibis.theo} 
is a reformulation of theorem~\ref{projectivesubform.theo}. 
Hence, we may consider it as an analogue for algebras with involution 
of the subform theorem up to similarities. 
It would be nice to have a proof which does not use any version of the 
subform theorem. 

As opposed to this, we do not consider theorem~\ref{scai.theo} as an analogue 
of the subform theorem. 
One should notice, in particular, that condition (iii) do imply that 
$(A,\sigma)$ contains a direct summand isomorphic 
to $(B,\tau)$, but this need not be $(eAe,\sigma_{|eAe})$ as the proof above shows. 
While under condition~(\ref{algebra.cond}) $(eAe,\sigma_{|eAe})$ 
is isomorphic to $(B,\tau)$, 
where $e$ precisely is the idempotent mentioned in the
condition. Hence, theorem~\ref{scai.theo} actually is a criterion of
isomorphism rather than a subform theorem. 
\end{rema} 

\section{Hermitian forms and the
  index $2$ case} 
\label{ind2.sect}

Before proving the theorem for algebras of index $2$, we translate
the conditions of theorem~\ref{scaibis.theo} in terms of hermitian
forms. 

\subsection{Conditions (ii) and (iii) in terms of
  hermitian forms}

\begin{prop} \label{hermitian.prop} 
Let $(A,\sigma)$ and $(B,\tau)$ be two Brauer equivalent 
central simple algebras both endowed with an orthogonal involution. 
Fix representations $(A,\sigma)=(\End_D(M),\ad_{h_M})$ and 
$(B,\tau)=(\End_D(N),\ad_{h_N})$
and let $q_V$ and $q_W$ 
denote as in \S\ref{scai.sec} the quadratic forms corresponding via Morita theory 
to the hermitian forms ${h_M}_{F_B}$ and ${h_N}_{F_B}$. 
We denote by $\delta$ a fixed non zero element in the generic line 
$L_B$. 

The algebras with involution  $(A,\sigma)$ and $(B,\tau)$ 
satisfy condition (iii) if and only if 
there exists an anisotropic vector $n\in N$,  
a vector $m\in M$ and a scalar $\la\in F^\times$ such that 
$h_M(m)=\la h_N(n)$, and $q_V$ represents $\la q_W(\delta)$. 
\end{prop} 

\begin{rema} 
Note that in particular, if $(A,\sigma)$ and $(B,\tau)$ satisfy 
condition (iii), then the hermitian forms $h_M$ and $h_N$ have 
a common value up to a central scalar $\la\in F^\times$. 
\end{rema} 

\begin{proof} 
Let us assume first that $(A,\sigma)$ and $(B,\tau)$ satisfy condition (iii). 
We use the computations made in 
the proof of theorem~\ref{scai.theo}, and keep the same notations. 
Consider $J$, $e$ and $\Psi$ as given by condition (iii). 
We denote again by $M_0$ the image of $e$ and we let 
$\psi:\,M_0\tilde\ra M$ be an isomorphism such that 
$\Psi(f)=\psi f\psi^{-1}$ for any $f\in A$. 
As recalled in remark~\ref{ideal.rema}, the ideal $J\subset B=\End_D(N)$ 
can be written as $J=\Hom_D(N/{(n.D)^\perp},N)$, for some 
anisotropic vector $n\in N$. 
Define $m$ to be $m=\psi^{-1}(n)$. 
Since $n$ is anisotropic, we may choose an orthogonal basis 
$n_1,\dots,n_2$ of $(N,h_N)$ over $D$ such that $n_1=n$. 
One may then easily check that the elements $g_{1,j}$ and $g_1$ 
defined in the proof of theorem~\ref{scai.theo} belongs to 
$J_{F_B}\cap I_B$, and condition (iii) proves the existence 
of corresponding elements $f_{1,j}$ and $f_1$ in $A_{F_B}$. 
By the computations made in the proof of theorem~\ref{scai.theo} 
we get that $h_M(m)=\la h_N(n)$ for some $\la\in F^\times$. 

Moreover, if we let $\nu=f_{1,1}(m_1\otimes 1\otimes e_1)\in V_0$, 
equation (\ref{q(nu).equa}) proves that 
$q_V(\nu)=b_{q_{V_0}}(f_{1,1}(m_1\otimes 1\otimes e_1))
=\la q_W(\delta)$. 

Let us now prove the converse. Since $n$ is anisotropic, 
there exists orthogonal basis $m_1,\dots m_r$ of $M$ and 
$n_1,\dots,n_s$ of $N$ such that $n_1=n$ and $m_1=m$. 
We let $J$ be $J=\Hom_D(N/{(n.D)^\perp},N)$, 
and take $e$ to be the orthogonal projection on the submodule 
$M_0$ of $M$ generated by $m_0,\dots, m_s$ and 
$\psi$ to be the isomorphism $\psi:\, M_0\tilde\ra N$ 
given by $\psi(m_i)=n_i$. 

Any endomorphism $g\in J_{F_B}\cap I_B$ is defined by 
\[
g(n_1\otimes 1\otimes e_j)=\alpha_j \delta\text{ for some }\alpha_j\in F_B
\text{ and }g(n_i\otimes 1  \otimes e_j)=0\text{ for }i\geq 2. 
\]
We then define $f\in A_{F_B}=\End_{F_B}(V)$ by 
\[
f(m_1\otimes 1\otimes e_j)={\alpha_j}\nu\\
\text{ and }f(m_i\otimes 1\otimes e_j)=0\text{ for }i\geq 2,
\] 
where 
$\nu\in V$ satisfies $q_V(\nu)=\la q_W(\delta)$. 
As in lemma~\ref{tau(g)g.lemm}, one may check that 
\[
\sig(f)f(m_1\otimes 1\otimes e_j)={\alpha_j}q_V(\nu)
(m_1h_M(m_1)^{-1}\otimes 1\otimes \sum_{l=1}^d {\alpha_l}\frac{e_l}{a_l}),
\] 
\[\text{ and }
\sig(f)f(m_i\otimes 1\otimes e_j)=0\text{ for }i\geq 2.
\]  
Hence we get 
\[\psi e\sig(f)fe\psi^{-1}(n_1\otimes 1\otimes e_j)
={\alpha_j}\la q_W(\delta)
(n_1\la^{-1}h_N(n_1)^{-1}\otimes 1\otimes \sum_{l=1}^d {\alpha_l}\frac{e_l}{a_l}),
\]
\[\text{and }\psi e\sig(f)fe\psi^{-1}(n_i\otimes 1\otimes e_j)=0\text{ for }i\geq 2
\] 
which by lemma~\ref{tau(g)g.lemm} implies 
$\psi e\sig(f)fe\psi^{-1}=\tau(g)g$ and hence finishes the proof. 
\end{proof} 

From the previous proposition 
and the correspondence between left ideals $J$ of reduced dimension $D$ 
such that $\tau(J)$ is anisotropic and anisotropic vectors in $(N, h_N)$ 
recalled in remark~\ref{ideal.rema}, we deduce a translation of
condition (ii) in terms of hermitian forms: 
\begin{coro} 
Let $(A,\sigma)$, $(B,\tau)$, $(M,h_M)$, $(N,h_N)$ and $\delta$ be as
in proposition~\ref{hermitian.prop}. 

The algebras with involution $(A,\sigma)$ and $(B,\tau)$ 
satisfy condition (ii) if and only if 
for any anisotropic vector $n\in N$, there exist  
a vector $m\in M$ and a scalar $\la\in F^\times$ such that 
$h_M(m)=\la h_N(n)$ and $q_V$ represents $\la q_W(\delta)$. 
\end{coro}

\subsection{Proof of theorem~\ref{scaibis.theo} in the index $2$ case} 
\label{index2.section} 

The theorem in the index $2$ case now follows
from the following proposition : 
\begin{prop} 
Let $(D,\gamma)$ be a division algebra of index at most $2$ 
with orthogonal involution 
and $(M,h_M)$ and $(N,h_N)$ two hermitian forms over $(D,\gamma)$, 
with $h_M$ anisotropic. 
Then $h_M$ contains $h_N$ as a subform if and only if 
$q_V$ represents $q_W(\delta)$. 
\end{prop} 

\begin{proof} 
If $D$ is split, this statement is the projective version of 
the classical subform theorem which 
is already known. Hence we assume $D$ is a division quaternion algebra. 
The condition is clearly necessary. 

To prove it also is sufficient, let us first extend scalars to $F_D$, so that 
the situation is split. 
By a result due both to I. Dejaiffe~\cite{Dej:01} and Parimala-Sridharan-Suresh~\cite{PSS}, 
we know 
that the involution $\sigma$ remains anisotropic after scalar extension 
to $F_D$, so that the quadratic form $q_V$ is anisotropic. 
Hence we can apply the subform theorem to $q_V$ and $q_W$. 
Since the generic point of a variety maps under scalar extension to the 
generic point of the extended variety, the condition implies that 
${q_V}_{F_D({\mathbb P}W)}$ represents ${q_W}^{proj}$, and hence 
$q_V$ contains $q_W$ as a subform. 
By Morita equivalence, this precisely means that 
$(h_M)_{F_D}$ contains $(h_N)_{F_D}$ as a subform, and 
it only remains to go down using the injectivity of 
the natural map from the Witt group of $D$ to the Witt group of 
$F_D$ and the excellence result of~\cite{PSS} 
\end{proof} 

{\bf Acknowledgements.} I would like to thank Philippe Gille, Bruno
Kahn, Parimala and Jean-Pierre Tignol for useful discussions on this
question while this work was in progress. 

\backmatter

\bibliographystyle{smfalpha} 
\bibliography{biblio-eng.bib} 

\end{document}